\magnification 1200
\documentstyle{amsppt}

\pagewidth{16.1 truecm}
\pageheight{22 truecm}
\hcorrection{-0.3 truecm}
\vcorrection{0.1 truecm}

\parskip 5pt
\document
\NoBlackBoxes
\nologo

\define\bo{\Bbb B^n}
\define\a{\alpha}

\define\dis{1-|\zeta|^2}
\define\diz{1-|z|^2}
\define\disak{1-|a_k|^2}
\define\disaj{1-|a_j|^2}

\define\akbaj{|1-\bar a_k a_j|}

\define\akbz{|1-\bar a_k z|}

\define\bak{\{v_k\}}
\define\sbk{\{v_k\}}
\define\sak{\{a_k\}}

\define\sumk{\sum_{k=1}^\infty}
\define\sumj{\sum_{j=1}^\infty}
\define\sumjk{\sum_{j:j\neq k}}

\define\bap{B_{\a}^{p}}

\define\ibo{\int_{\Bbb B^n}}

\define\nap{\frac{n+1}p+\alpha}
\define\lap{\ell_{\frac{n+1}p+\alpha}^p}
\redefine\phi{\varphi}

\define\bapp{B_{\a'}^{p'}}
\define\isf{\int_S}
\define\napp{\frac{n+1}{p'}+\alpha'}
\def\laprime{\ell_{\frac{n+1}{p'}+\alpha'}^{p'}}
\def\naprime{\frac{n+1}{p'}+\alpha'}
\def\uap{\frac{1}p+\alpha}
\def\uaprime{\frac{1}{p'}+\alpha'}
\define\naf{\frac{n+1}p+\a}
\define\nafp{\frac{n+1}{p'}+\alpha '}
\define\pp{\frac{p'p}{p-p'}}
\define\mlk{|\lambda_k|}

\topmatter
\title
Interpolating sequences for weighted Bergman
spaces of the ball
\endtitle

\author
Miroljub Jevti\'c, Xavier Massaneda, Pascal J. Thomas
\endauthor

\address
Miroljub Jevti\'c:
Matemati\v cki fakultet, Studentski trg 16,
11000 BEOGRAD, Yugoslavia
\endaddress
\email XPMFM39\@yubgss21.rzs.bg.ac.yu \endemail

\address
Xavier Massaneda: 
Departament de Matem\`atiques,
Universitat Aut\`onoma de Barcelona,
08193 BELLATERRA, Spain 
\endaddress
\email
xavier\@manwe.mat.uab.es
\endemail

\address
Pascal J. Thomas:
UFR MIG,
Universit\'e Paul Sabatier
118 route de Narbonne,
31062 TOULOUSE CEDEX,
France
\endaddress
\email
pthomas\@cict.fr
\endemail

\abstract Interpolating sequences for weighted Bergman spaces $\bap$,
$0<p\leq\infty$, $\a\geq-1/p$ are studied. We show that the natural
inclusions between $\bap$ for various $p$ and $\a$ are also verified
by the corresponding spaces of interpolating sequences. We also give
conditions (necessary or sufficient) for the $\bap$-interpolating
sequences. These are similar to the known conditions for the spaces
$H^p$ and $A^{-\a}$, which in our notation correspond respectively to
the particular cases $\a=-1/p$ and $p=\infty$.
\endabstract

\thanks
Second author partially supported by DGICYT grant PB92-0804-C02-02
\endthanks

\endtopmatter

\centerline{\it \S 0. Introduction.}

Let $\bap$ be the space of $f$ holomorphic in the
 unit ball $\bo$ of $\Bbb C^n$ such that
$(\diz)^\a f(z) \in L^p(\bo)$, where $0<p\leq\infty$, $\a\geq-1/p$
(weighted Bergman space). 
In this paper we study 
the interpolating sequences for various $\bap$. 
The limiting
cases $\a=-1/p$ and $p=\infty$ are respectively the
Hardy spaces $H^p$ and $A^{-\a}$, the holomorphic 
functions with polynomial growth of
order $\a$, which have generated particular
interest. Note that the class of spaces we are considering
is invariant under restriction to balls of lower complex
dimension, which justifies the choice of those
special weights.

As far as we know, for $n>1$ 
the first research on this subject was carried out by Amar for the
classical Bergman spaces, which in our notation correspond to the
case $\a=0$ (\cite{Am}). Amar's main result states that
separated sequences (in terms of the Gleason invariant distance) 
can be written as a finite union of interpolating sequences for $B_0^p$. 

A sufficient condition due to Berndtsson is known
for the case $H^\infty$ (\cite{Be}). Also, Varopoulos showed that if
$\sak_k$ is $H^\infty$-interpolating then $\sum_k (\disak)^n
\delta_{a_k}$ is a Carleson measure (\cite{Va}). Later 
the third author proved
that the same necessary condition holds for $H^1$, and it actually 
characterizes the finite unions of $H^1$-interpolating sequences
(\cite{Th1}). 

On the other hand, after Seip's characterization of
$A^{-\a}$-interpolating sequences in the unit disc (\cite{Se},
see another proof in \cite{Be-Or}), 
the second author
obtained some results for the case $n>1$ (\cite{Ma}). In particular,
$\sak_k$ is a finite union of $A^{-\a}$-interpolating sequences if
and only if $\sum_k (\disak)^{n+1}$ is a $(n+1)$-Carleson measure, or
equivalently, if and only if $\sak_k$ is a finite union of separated sequences.

It is worth noting that in \cite{Se}, Seip also 
implicitly gives a characterization
of interpolating sequences for all weighted Bergman spaces in the
disk. We spell out the details for the reader's convenience in an
appendix (\S 5).

Here we deal with different aspects concerning $\bap$-interpolating
sequences. In \S 1 we first collect some definitions and well-known
facts about weighted Bergman spaces and then introduce the natural
interpolation problem, along with some basic properties. In \S 2 we
describe in terms of $\a$ and $p$ the inclusions between $\bap$
spaces, and in \S 3 we show that most of these inclusions also hold for
the corresponding spaces of interpolating sequences. Unfortunately
our proof does not capture the intuitive conjecture given in
\cite{Th1} to the effect that for $p'<p$ every $H^p$-interpolating sequence
is also $H^{p'}$-interpolating. \S 4 is devoted to sufficient
conditions for a sequence to be $\bap$-interpolating,
expressed in the same terms as the conditions given in
\cite{Th1} for the Hardy spaces and in \cite{Ma} for $A^{-\a}$. In
particular we show, under some restrictions on $\a$ and $p$, that
finite unions of $\bap$-interpolating sequences coincide with finite
unions of separated sequences.

\ {}

\centerline{{\it \S 1. Definitions and first properties.}}

{\bf 1.1. Notations.}
For $z$, $w \in \Bbb C^n$, we set $z \bar w := \sum_{j=1}^n z_j {\bar w}_j$,
$|z|^2 := z \bar z$, and the unit ball
$\bo := \{ z \in \Bbb C^n : |z| < 1 \}$.

Given $a \in \bo$, $\varphi_a$ is
the involutive automorphism of the ball exchanging $0$ and $a$
(see \cite{Ru, 2.2.2}).
For $a$, $b \in \bo$, $d(a,b) := |\varphi_a (b)| =  |\varphi_b (a)|$
is the invariant distance between $a$ and $b$. Recall that
$$
1-d (a,b)^2 = 1- |\varphi_a (b)| ^2 =
\frac{(1-|a|^2)(1-|b|^2)}{|1- a \bar b |^2} .
$$
We call 
{\it hyperbolic balls} the sets
$E(z,r):=\{\zeta\in\bo : d(z,\zeta)<r\}$. 

The normalized Lebesgue measures on the ball and the sphere will
be denoted by $d m $ and $d \sigma$ respectively 
($d m_{2n}$ and $d \sigma_{2n-1}$ when we want to
stress the dimension).
The  measure
$d\tau (z):=(\diz)^{-(n+1)} dm(z)$ is invariant under
the automorphisms of the ball \cite{Ru, 2.2.6}.
In particular, $\tau (E(z,r))$ depends only on $r$.

Throughout this paper we will be using the following estimates:

\proclaim{\smc Lemma 1.1} Let $a,b\in\bo$, $c>0$, $t>-1$. Then
\roster
\item"(a)" $\dsize{\ibo\frac{(\diz)^t}{|1-\bar a z|^{n+1+c+t}} \, dm(z)}
\simeq (1-|a|^2)^{-c}$.

\item"(b)" $\dsize{\int_S\frac{d\sigma(\zeta)}{|1-\bar a\zeta|^{n+c}} \, dm(z)}
\simeq(1-|a|^2)^{-c}$.
\item"(c)" $\dsize{
\ibo \frac{(1-|z|^2)^t}{|1-z\bar a|^{n+1+c+t}|1-z\bar b|^{n+1+c+t}} \, dm(z)
}$
\newline
$\preceq |1-a\bar b|^{-(n+1+c+t)} \{ \min( 1-|a|^2, 1-|b|^2 ) \}^{-c}$.
\endroster
\endproclaim

\demo{\sl Proof} (a) and (b) are given in \cite{Ru, 1.4.10}. To prove (c)
split into the cases $|1-z\bar a| \geq \frac{\sqrt2}2 |1-a\bar b|$
and $|1-z\bar a| \leq \frac{\sqrt2}2 |1-a\bar b|$, which implies
$|1-z\bar b| \geq \frac{\sqrt2}2 |1-a\bar b|$ by the triangle
inequality \cite{Ru, 5.1.2(i)}; then apply (a). \qed
\enddemo

{\bf 1.2. Weighted Bergman Spaces.}
For $p>0$, $\a \in \Bbb R$,
let $L_\a^p (\bo)$ be the space of all measurable
complex-valued functions on $\bo$ such that 
$(1-|z|^2)^{\a } f(z) \in L^p$, i.e. for $p<\infty$,
$$
\| f \|_{p,\a}^p := \ibo (1-|z|^2)^{\a p} |f(z)|^p \, dm(z) < \infty ,
$$
and for $p=\infty$,
$$
\| f \|_{\infty,\a} := \sup_{z \in \bo} (1-|z|^2)^{\a } |f(z)| < \infty .
$$
For $p>0$, $\a > -1/p$, and 
denoting by $H(\bo)$ the space of
holomorphic functions in the ball, 
the {\it weighted Bergman space} is $\bap (\bo) := H(\bo) \cap L_\a^p (\bo)$.

For $\a \leq -1/p$, the above condition 
only holds for the zero function, but we sometimes will
use the limiting case of the {\it Hardy spaces}:
$$
B_{-\frac1p}^p := H^p (\bo) = \{ f \in H(\bo) : 
\| f \|_{H^p} = \| f \|_{p,-\frac1p}^p := 
\sup_{0<r<1} \int_{\partial \bo} |f(r \zeta)|^p \, d \sigma (\zeta) 
< \infty \}.
$$
The facts presented in this section are
essentially well-known, but we recap them here for the
reader's convenience, and also to write them in our
notation, which differs from 
those of Horowitz \cite{Ho}, Coifman and Rochberg \cite{Co-Ro},
\cite{Ro}, and Seip \cite{Se}. 

The following statement is an immediate consequence of \cite{Ru, pag. 14}.

\proclaim{\smc Lemma 1.2} 
If $\ell \in \Bbb Z_+$ and
$f \in B_{\a - \frac{\ell}p}^p ( \Bbb B^{n+\ell}) $,
then its restriction to $\bo \times \{ 0 \}$ lies in
$\bap (\bo)$; conversely whenever $g \in \bap (\bo)$,
its trivial extension (constant along the vertical directions)
must be in $B_{\a - \frac{\ell}p}^p ( \Bbb B^{n+\ell}) $. 
\endproclaim


\proclaim{\smc Lemma 1.3}
For $p>0$, $\a \geq -1/p$, there exists a constant
$c = c(\a, p, n) >0$ such that for all $z \in \bo$
$$
|f(z)| \leq c \| f \|_{p,\a} (1-|z|^2)^{-(\nap)},
$$
and this is the best possible exponent.
\endproclaim

\demo{\sl Proof}
Lemma 1.2 allows us to reduce this to the case $n=1$
by considering the disk through $0$ and $z$. Then use
the mean value inequality on the disk $D(z, \frac12 (1-|z|))$.

That the estimate is sharp can be seen by considering
 the functions
$f_{N,a} (z) := (1-z \bar a )^{-N}$ for $a \in \bo$,
$N> \nap$. Lemma 1.1 shows that 
$\| f_{N,a} \|_{p,\a} \approx (1-|a|^2)^{\nap - N}$. \qed
\enddemo

Lemma 1.3 says that $\bap \subset B_{\nap}^\infty$.
A more complete catalogue of inclusions will be given in Section 2.

\proclaim{\smc Lemma 1.4}
There exists $c' = c'(\a, p, n) > 0$ such that for any
$a$, $b \in \bo$ with $d(a,b) < 1/2$,
$$
| f(a) - f(b) | \leq c' \| f \|_{p,\a} (1-|a|^2)^{-(\nap)} d(a,b) .
$$
\endproclaim

\demo{\sl Proof}
We apply the generalized Schwarz Lemma (see \cite{Ru, 8.1.4})
over $E(a, \frac12 )$ together with the estimate from Lemma 1.3. \qed
\enddemo

{\bf 1.3. The interpolation problem.}
We'd like to call a sequence $\sak_k \subset \bo$ 
{\it interpolating} for the space $\bap$ if, given
arbitrary values $\bak_k$ subject to some reasonable restrictions,
there exists $f \in \bap$ such that $f(a_k ) = v_k$ for
all $k \in \Bbb Z_+$. It is not so easy in general to
determine what is a reasonable set of possible values
for the restrictions of holomorphic functions to a sequence
of points, see \cite{Br-Ni-\O y} for the case of $H^p (\Bbb D)$.

Here we want to impose  growth restrictions only
on the sequence $\bak_k$, meaning that if 
$\bak_k$ is in the set of possible values 
(required to be a Fr\'echet space),
then so is $\{v'_k\}_k$ whenever $|v'_k| \leq |v_k| $ for all
$k$. Since we'd like to allow the values 
$v_k^N := (1-|a_j|^2)^{-(\nap) + N} f_{N,a_j} (a_k)$, for any
given $j$, and since
$ \| (1-|a_j|^2)^{-(\nap) + N} f_{N,a_j} \|_{p,\a} $ does not depend
on $j$, 
we shall adopt the following provisional
definition:

\proclaim{\it Definition} 
We say that $\sak_k$ has the property (*) for $\bap$ iff there
exists $M>0$ such that for any $j \in \Bbb Z_+$,
there exists $f_j \in \bap$ such that $\|f_j\|_{p,\a} \leq M$
and $f_j (a_k) = 0$ for $k \neq j$ and 
$f_j (a_j) = (1-|a_j|^2)^{-(\nap)}$.
\endproclaim

The requirement on the norm could be deduced from our more
general requirements if we were to assume that the 
sequences $\{ v_k^N \}_k$ lie in the unit ball of some
Banach space, which would be reasonable in the case
where $p \geq 1$.

At any rate, (*) in itself leads to stronger properties,
starting with a geometric restriction on the sequence $\sak_k$.

\proclaim{\it Definition}
For $p>0$, $\beta \in \Bbb R$, let
$$
\ell_\beta^p =
\ell_\beta^p (\sak) := \bigl\{ \bak_k \subset \Bbb C \ :\ 
\{ (1-|a_k|^2)^\beta v_k \}_k \in \ell^p \bigr\} 
$$
and $\| v \|_{p,\beta}^p := \sum_k [(1-|a_k|^2)^\beta |v_k|]^p$,
$\; \| v \|_{\infty,\beta} := \sup_k [(1-|a_k|^2)^\beta |v_k|]$.
\endproclaim

\proclaim{\it Definition}
We say that a sequence $\sak_k$ is {\it separated} iff there
exists $\delta > 0$ such that for any $j \neq k$,
$d(a_j , a_k) \geq \delta$.
\endproclaim

\proclaim{\smc Lemma 1.5}
For any $p>0$, $\a > -1/p$, 
and any separated sequence $\sak_k$, the restriction map
$f \mapsto \{f(a_k)\}_{k\in \Bbb Z_+}$ is bounded
from $\bap$ to $\lap$.
\endproclaim

\demo{\sl Proof}
For $p= \infty$, this is trivial.

For $p<\infty$, notice that the separatedness implies that
for some $\delta >0$
the hyperbolic balls
$E(a_k, \delta)$ are pairwise disjoint. Thus
$$
\| f \|_{p, \a}^p \geq
\sum_k \int_{E(a_k, \delta)} (1-|z|^2)^{\a p} |f(z)|^p \, dm(z)
\gtrsim \delta^{n+1} \sum_k (1-|a_k|^2)^{\a p + n+1} |f(a_k)|^p ,
$$
using the plurisubharmonicity of $|f|^p$. \qed
\enddemo

Assume conversely that $p > 0$, $\a \geq -1/p$ and that $\sak_k$
verifies (*) for $\bap$. Applying Lemma 1.4 to the $f_j$ given by
(*), we find that
$$
(1-|a_j|^2)^{-(\nap)} \leq c' M (1-|a_j|^2)^{-(\nap)} d(a_j,a_k)\ ,
$$
so $\sak_k$ is separated.

This motivates the following definition:

\proclaim{\it Definition}
We say that $\sak_k$ is an {\it interpolating sequence} for
$\bap$, denoted by $\sak \in Int (\bap)$, iff 
for any $\bak \in \lap$, there exists $f\in \bap$ 
such that $f(a_k) = v_k$ for all $k$.
\endproclaim

This definition appeared in \cite{Sh-Sh} in the case
$\a = -1/p$, $n=1$, and occurs under various guises in
\cite{Am}, \cite{Ro}, \cite{Se} and \cite{Th1}, \cite{Th2}.

Applying Baire's Theorem to the closure of the images
under the restriction map of balls of arbitrarily large radius,
we see that whenever $\sak_k$ is an interpolating sequence,
there is an $M>0$ such that 
for any $\bak \in \lap$, for any $\varepsilon > 0$, 
there exists $f\in \bap$, $\| f \|_{p,\a} \leq M$,
such that $\| f(a_k) - v_k \|_{p, \nap} < \varepsilon $.
Thus, applying again Lemma 1.4, we have that $\sak_k$ is separated,
hence for $\a > -1/p$ the restriction map is bounded. 
For $\a = -1/p$, $p\geq 1$, the proof of \cite{Th1, Theorem 2.2} shows that 
(*) for $H^p$ implies that $\sum_k (\disak)^n \delta_{a_k}$ is a 
Carleson measure (see \cite{Th2} for a direct proof when $p>1$). 
Applying Baire's Theorem likewise, we 
can correct \cite{Th1} to get the boundedness 
of the restriction mapping claimed there.

We can then apply the Open Mapping Theorem.

\proclaim{\smc Corollary 1.6}
If $\sak_k$ is an interpolating sequence for
$\bap$, then $\sak_k$ is separated, the restriction map is 
bounded from $\bap$ to $\lap$, and 
there exists a constant $M>0$ so that the function $f$ in
the definition can be chosen with the additional condition
$\| f \|_{p, \a} \leq M \| v \|_{p, \nap}$.
\endproclaim

The constant $M$ is called {\it constant of interpolation of $\sak_k$}.

{\bf 1.4. Invariance under automorphisms and restriction to subspaces.}
For $\varphi$ any automorphism (holomorphic self-map)
of the ball, let
$$
T_\varphi f (z) := \left( 
\frac{1-|\phi^{-1} (0) |^2}{(1-z \overline{\phi^{-1} (0)})^2} \right)^{\nap}
f \circ \phi (z).
$$

\proclaim{\smc Lemma 1.7}
\roster
\item"(a)" $T_\phi$ is an isometry of $\bap$.
\item "(b)" If $\varphi$ is an automorphism
of the ball, and 
$\sak_k$ is an interpolating sequence for $\bap$, 
then so is $\{ \phi (a_k ) \}_k$, with the same constant
of interpolation.
\endroster
\endproclaim

\demo{\sl Proof}
(a) is trivial when $p=\infty$.
Any automorphism is a composition of a map $\phi_a$ and
a rotation, and the result is immediate in the latter case.
For the former, in the case where $p < \infty$, $\a > -1/p$,
$$
\align
\ibo |T_{\varphi_a} f (z) |^p (1-|z|^2)^{\a p} \, dm(z) 
&=
\ibo (1-|\varphi_a (z) |^2 )^{n+1+\a p}
|f (\phi_a(z)) |^p \, d\tau(z) \\
&=
\ibo (\dis)^{n+1+\a p} |f (\zeta) |^p \, d \tau (\zeta)
= \| f \|_{p,\a}^p,
\endalign
$$
which finishes the proof since $(T_\phi)^{-1}= T_{\phi^{-1}}$,
as can be seen by an elementary calculation using 
\cite{Ru, 2.2.5} and the fact that $|\phi^{-1} (0) |=|\phi (0) |$. 

(b) Take $v \in \ell_\beta^p (\{ \phi (a_k ) \})$. Then 
$$ 
\left\{ \left( 
\frac{1-|\phi^{-1} (0) |^2}
{(1-a_k \overline{\phi^{-1} (0)})^2} \right)^{\nap} 
v_k \right\}
\in \ell_\beta^p (\sak) ,
$$
with the same norm,
so there is $F$ such that 
$\| F \|_{p, \a} \leq M \| v \|_{p, \nap}$ and
$$
F (a_k ) = \left( 
\frac{1-|\phi^{-1} (0) |^2}
{(1-a_k \overline{\phi^{-1} (0)})^2} \right)^{\nap} 
v_k \ .
$$
Then $G:= T_{\phi^{-1}} F$ solves the original problem,
with $\| G \|_{p, \a} \leq M \| v \|_{p, \nap}$. \qed
\enddemo

The following lemma, which follows immediately from Lemma 1.2, 
has been used in \cite{Am}, and provides some necessary conditions 
for a sequence to be interpolating.

\proclaim{\smc Lemma 1.8}
Suppose $\sak \subset \bo \times \{ 0 \} \subset \Bbb B^{n+\ell}$,
where $\ell \in \Bbb Z_+^*$. Let $\a \geq (\ell-1)/p$.
Then $\sak \in Int (\bap (\bo))$ if and only if
$\sak \in Int (B_{\a - \frac{\ell}p}^p ( \Bbb B^{n+\ell}) )$.
\endproclaim

\par
\ {}
\par

{\bf 1.5. Stability.} 
The proof of the following Lemma was sketched in \cite{Lu, \S 6.II}.

\proclaim{\smc Lemma 1.9}
For $0< p \leq \infty$ and $\a> -1/p$
or $1\leq p \leq \infty$ and $\a \geq -1/p$,
let $\sak\in Int(\bap) $ and let $\{a_k'\}_k $ be another sequence in 
$\bo$. There exists $\delta>0$ such that if
$$
d(a_k, a_k')<\delta\qquad\forall k\in\Bbb Z_+
$$
then $ \{a_k'\}\in Int(\bap) $.
\endproclaim

\demo{\sl Proof} 
Case $\a > -1/p$.
Let $v\in \ell_{\frac{n+1}p+\a}$. Denote $a=\sak_k$,
$a'=\{a_k'\}_k$ and $v^0=v $. By hypothesis there exists $f_0\in\bap$ such
that $f_0(a)=v^0$ ($f_0(a_k)=v_k^0$ for all $k$) and $\|f_0\|_{p,\a}\leq
\|v^0\|$, where $M$ denotes the constant of interpolation  of $\sak_k$.
Consider now $v^1:=v^0-f_0(a')$. 

{\sl Claim}. For $\delta$ small enough,
$\|v^1\|\leq \gamma \|v^0\|$, with  $\gamma<1$. 

To see this we use a general estimate for holomorphic functions
which is a refinement of Lemma 1.4. Let $f$ be
holomorphic and let $z,w\in\bo$ with $d(z,w)<r <1$. The
plurisubharmonicity of $|f(z)-f(w)|^p$ as a function of $z$ together with
a gradient estimate shows that there exists a constant $C= C(r)>0$ such that
(see \cite{Lu, Lemma 3.1} or \cite{Th2, Lemma 2.4.4}): 
$$
|f(z)-f(w)|^p\leq C d^p(z,w) \int_{E(w,r)} |f(\zeta)|^p\; d\tau(\zeta)
$$
for any $r>\frac 23 d(z,w)$. With this estimate applied to $f_0$ we have,
provided that $r$ is chosen small enough so that the
invariants balls $E(a_k,r)$ are pairwise disjoint:
$$
\multline
\|v^1\|^p=\sum_k(\disak)^{n+1+\a p}|f_0(a_k)-f_0(a_k')|^p\leq\\
\leq C \sum_k(\disak)^{n+1+\a p} d^p
(a_k,a_k')\int_{E(a_k,r)}|f_0(\zeta)|^p\; d\tau(\zeta)\leq\\
\leq C\delta^p \sum_k \int_{E(a_k,r)}|f_0(\zeta)|^p (\dis)^{\a p}\;
dm(\zeta)\leq C\delta^p \|f_0\|_{p,\a}^p\leq C\delta^p M^p \|v^0\|^p\ . 
\endmultline
$$
Choosing $\delta$ so that $\gamma^p:= C\delta^p M^p<1$ the claim is
proved. 

Take now $f_1\in\bap$ with $f_1(a)=v^1$ and $\|f_1\|_{p,\a}\leq M\|v^1\|$,
and define $v^2=v^1-f_1(a')$. An iteration of this construction provides
functions $f_j\in\bap$ with $f_j(a)=v^j=v^{j-1}-f_{j-1}(a')$ and
$\|f_j\|_{p,\a}\leq M \|v^j\|\leq M\gamma^ j \|v^0\|$. Finally
the function $f=\sum_j f_j$ solves the interpolation
problem for $\{ a'_k \}_k$.

Case $\a = -1/p$, $p\geq1$.
We can use Luecking's estimate, together with the fact
that, when $\sak \in Int (H^p)$, $p\geq1$, the measure 
$\sum _k (\disak)^n \delta_{a_k}$ is a Carleson measure \cite{Th1,
Theorem 2.2}, and therefore so is the measure 
$\sum _k (\disak)^{-1} \chi_{E(a_k,r)} dm$. Then, applying Lemma 3.1
(see section \S 3) we get
$\sum _k (\disak)^n |f(a_k)-f(a_k')|^p \lesssim  \| f\|_{H^p}^p$
for any function for which the right hand side is finite.
The proof then proceeds as before. \qed

\enddemo

The results described in this section also hold for the corresponding
interpolating sequences for the spaces $b_\a^p$ of $\Cal M$-harmonic
(instead of holomorphic) functions in $L_\a^p$. This is so because
the main ingredients used above, namely Lemma 1.3 and the separatedness of the
interpolating sequences, can be proven likewise in the $\Cal
M$-harmonic case.

\par
\ {}
\par

\centerline{\it \S 2. Inclusions for $\bap$ and $\ell^p_\beta$.}

Our purpose is now
to describe in terms 
of the values $\a$ and $p$ the relationship between $\bap$ spaces.

\proclaim{\smc Lemma 2.1}
\roster
\item "(a)"
If $p\leq p'$, then $B_{\a}^{p}\subset\bapp$ iff 
$\ \a+\frac{n+1}{p} \leq \a'+\frac{n+1}{p'}\ $.
\item "(b)" 
If $\ p\geq p'\ $, then 
$\ \a+\frac 1{p}<\a'+\frac 1{p'}\ $
$\Rightarrow B_{\a}^{p}\subset\bapp$ 
$\Rightarrow \ \a+\frac 1{p} \leq \a'+\frac 1{p'}\ $.
\endroster

\endproclaim

In particular, if $\a  < \a' + 
\min \left( (n+1)(\frac{1}{p'} -\frac{1}{p}) , 
(\frac{1}{p'} -\frac{1}{p}) \right)$, then $B_{\a}^{p}\subset\bapp$,
and if $B_{\a}^{p}\subset\bapp$, then 
$\a \leq \a'+ 
\min \left( (n+1)(\frac{1}{p'} -\frac{1}{p}) , 
(\frac{1}{p'} -\frac{1}{p}) \right)$. Those results
have been obtained in the case $n=1$ by Horowitz \cite{Ho}. 

Among many other results,
Coifman and Rochberg 
\cite{Co-Ro, Propositions 4.2 and 4.4}
 prove that if
$p \leq p' \leq 1$,  $\a'$, $ \a \geq - \frac1{p}$
and $\a+\frac{n+1}{p} = \a'+\frac{n+1}{p'}$, then
$\bap \subset \bapp$. Their
proof is valid for a whole class of symmetric domains.

\demo{\sl Proof}
(a). Assume $\ \a+\frac{n+1}{p} \leq \a'+\frac{n+1}{p'}\ $
and  $f \in  B_{\a}^{p}$. The case $p' = \infty$ was
settled by Lemma 1.3. For $p'<\infty$, one has
$$
\ibo |f(z)|^{p'} (\diz)^{\a' p'} dm(z)\leq  c^{p'-p}\ibo|f(z)|^{p}
(\diz)^{(\frac{n+1}{p}+\a)(p-p')+\a' p'} dm(z)\ .
$$
This integral is controlled by $\Vert f\Vert_{\a,p}^{p}$ 
whenever $(\frac{n+1}{p}+\a)(p-p')+\a' p'\geq \a p$, 
that is, when $\a+\frac{n+1}{p}\leq\a'+\frac{n+1}{p'}$. 

Conversely, assume  $\a+\frac{n+1}{p} > \a'+\frac{n+1}{p'}$.
As in Section 1.2, let  
$f_{\gamma ,a} (z) = (1 - z \cdot \bar a ) ^{-\gamma}$.
Whenever 
$\gamma > \a+\frac{n+1}{p'}$, 
$\| f_{\gamma ,a}  \|_{\a',p'} 
\approx (1 - |a|^2 ) ^{\a'+\frac{n+1}{p'} - \gamma}$.
Choosing $\gamma > \a+\frac{n+1}{p}$, we see
that $\| f_{\gamma ,a}  \|_{\a',p'} / \| f_{\gamma ,a} 
\|_{\a,p} $ cannot be bounded as $|a|$ tends to $1$.

 (b). Assume $\ \a+\frac 1{p}<\a'+\frac 1{p'}\ $ and $f \in B_{\a}^{p}$. 
H\"older's inequality with 
exponent $\frac{p}{p'}\geq 1$ yields
$$
\multline
\ibo |f(z)|^{p'} (\diz)^{\a' p'} dm(z)\leq\\
\leq\left(\ibo|f(z)|^{p} (\diz)^{\a p} dm(z)\right)^{\frac 
{p'}{p}}\left(\ibo(\diz)^{(\a'-\a)\frac{p'p}{p-p'}} dm(z)\right)^{\frac{p-
p'}{p}}\ ,
\endmultline
$$
and by hypothesis both integrals are finite.
If $p'$ or $p$ is infinite, the analogous proof goes through
even more easily.

When $\a+1/p=0$, the hypothesis
becomes $f \in H^{p}$. Using 
$\isf |f|^{p'} d \sigma \leq \isf |f|^{p} d \sigma$
and integration in polar coordinates, we see that
$$
\ibo |f(z)|^{p'} (\diz)^{\a' p'} dm(z) \preceq 
\| f\|_{H^{p}}^p \int_0^1 (1-r^2 )^{\a' p'} d r < \infty ,
$$
since $\a' p' > -1$.

When $\a + 1/p = \a' + 1/p' >0$,
$p>p'$, it is possible to find
$f \in L_{\a}^{p} (\bo) \setminus L_{\a'}^{p'} (\bo)$. 
However we don't know whether there exist holomorphic functions
with that property. 

When $\a + 1/p > \a' + 1/p'$, we will construct
an example in the following way. For $\kappa > 1$
and for $0<r<1$, 
choose a set $\{ \eta_k \}_k \subset S$ maximal for the property
that the Koranyi balls 
$$
K(\eta_k , \kappa r ) :=
\{ \zeta \in S : |1-\zeta \overline{\eta_k}| < \kappa r \}
$$
are disjoint. Recall that the quantity $ |1-\zeta \bar \xi|^{1/2} $
verifies the triangle inequality over $\overline{\bo}$
\cite{Ru, 5.1.2}. For $\gamma>0$, let
$$
F_{\gamma,r}(z):= \sum_k \frac1{(1-(1-r)z \overline{\eta_k})^\gamma} .
$$
{\it Claim:}
 For $\kappa$ large enough,
and for all $0< p' \leq \infty$, $\a' \geq -1/p'$ and 
$\gamma > \max(\napp, n)$, then
$\| F_{\gamma,r} \|_{p',\a'} \approx r^{\frac1{p'} + \a' - \gamma}$,
where the constants involved may depend on $p'$, $\a'$,
$\gamma$, but not on $r$.

>From the claim we conclude as above,
letting $r \to 0$, that $B_{\a}^{p}\not\subset\bapp$
when $\a + 1/p > \a' + 1/p'$.

{\it Proof of Claim.} We need to perform a pointwise estimate on 
$F_{\gamma,r}(z)$. 
Suppose $z=\rho \zeta$, $\zeta \in S$. Suppose first
that $\zeta \in K(\eta_k , r )$. Then for any $j\neq k$,
$$
\inf_{\xi \in K(\eta_j , \kappa r )} 
| 1 - (1-r) \rho \zeta \bar \xi | \geq 
C | 1 - (1-r) \rho \zeta \overline{\eta_j}| ,
$$
and using Riemann sums and Lemma 1.1,
$$
\align
|F_{\gamma,r}( \rho \zeta )| &\geq
\frac1{|1-(1-r) \rho \zeta  \overline{\eta_k}|^\gamma} 
- C (\kappa r)^{-n} \int_{S \setminus K(\eta_k , \kappa r )} 
\frac{d \sigma (\xi)}{(| 1 - (1-r) \rho \zeta \bar \xi |^\gamma}\geq\\
&\geq C \left( \max ( (1-\rho),r ) \right)^{-\gamma} 
- C \kappa^{-n} r^{-n} \left( \max ( (1-\rho),r ) \right)^{n-\gamma} ,
\endalign
$$
The last quantity is
bounded below by $C_1 r^{-\gamma}$ for $\kappa$ large enough
and $\rho \geq 1-r$. The estimate from below for $p=\infty$
is thus secured.

For $p'<\infty$ and $\a' = -1/p'$, we have that
$$
\sup_{\rho<1} \left( 
\int_S |f(\rho \zeta)|^{p'} d \sigma (\zeta) 
\right)^{\frac1{p'}} 
\gtrsim r^{-\gamma} \ ,
$$
and for $p'<\infty$ and $\a' > -1/p'$, 
$$
\ibo |F_{\gamma,r}(z)|^{p'} (\diz)^{\a' p'} \, dm(z)
\geq \int_{1-r}^1 r^{-\gamma p'} (1-\rho^2)^{\a' p'} \rho^{2n-1} \, d \rho
\approx r^{1+(\a' -\gamma)p'} .
$$
On the other hand, for any $\zeta \in S$,
by the same arguments,
$$\align
|F_{\gamma,r}( \rho \zeta )| &\leq
\sum_{k: \eta_k \in K(\zeta, 2 \kappa r)}
\frac1{|1-(1-r) \rho \zeta  \overline{\eta_k}|^\gamma} 
+ C (\kappa r)^{-n} \int_{S \setminus K(\zeta , \kappa r )} 
\frac{d \sigma (\xi)}{(| 1 - (1-r) \rho \zeta \bar \xi |^\gamma}
\\
& \leq C(n) \left( \max ( (1-\rho),r ) \right)^{-\gamma} 
+ C \kappa^{-n} r^{-n} \left( \max ( (1-\rho),r ) \right)^{n-\gamma} .
\endalign
$$
So for $\rho \geq 1-r$, $|F_{\gamma,r}(\rho \zeta)| \leq C_2 r^{-\gamma}$,
and for $\rho \leq 1-r$, 
$|F_{\gamma,r}(\rho \zeta)| \leq C r^{-n}(1-\rho)^{n-\gamma}$.

In the case $\a' > -1/p'$, assuming $\gamma>\frac{n+1}{p'}+\a'$, we
integrate this to get:
$$
\multline
\ibo  |F_{\gamma,r}(z)|^{p'} (\diz)^{\a' p'} \, dm(z) \preceq
\\
 \preceq 
\int_0^{1-r} r^{-np'} (1-\rho^2)^{(n-\gamma + \a')p'}\rho^{2n-1} \, d \rho +
\int_{1-r}^1 r^{-\gamma p'} (1-\rho^2)^{\a' p'}\rho^{2n-1} \, d \rho
\preceq\\
\preceq  r^{1+ (\a' -\gamma )p'}\ . \qed
\endmultline
$$
\enddemo

The spaces of possible values we are considering verify
similar inclusions.

\proclaim{\smc Lemma 2.2} Suppose $\sak_k$ is separated. Then
\roster
\item "(a)" If $p\leq p'$, then $\lap \subset \laprime$ if and only if
$\nap \leq \naprime$;

\item "(b)" If $p \geq p'$ and $\uap < \uaprime$, then $\lap \subset \laprime$.
\endroster
Conversely to (b),
\roster
\item "(c)" If $p >p'$ and $\uap \geq \uaprime$, then
there exists a separated sequence $\sak_k$ such that 
$\lap \not\subset \laprime$.
\endroster
\endproclaim

\demo{\sl Proof} (a) Left to the reader (similar to Lemma 2.1).

(b) First observe that the separatedness of the sequence
implies that for any $\varepsilon > 0$,
$\sum_k (\disak)^{n+\varepsilon} < \infty$ (see Lemma 4.1).
For $p < \infty$, applying H\"older's inequality,
$$
\multline
\sum_k (\disak)^{n+1 + \a' p'} | v_k |^{p'} 
\leq 
\left( 
\sum_k 
\left( (\disak)^{\a p'} |v_k|^{p'} 
\right)^{\frac{p}{p'}} (\disak)^{n+1}
\right)^{\frac{p'}{p}}
\times\\
\times
\left( 
\sum_k 
\left( (\disak)^{(\a' - \a)p'}
\right)^{\frac1{1-\frac{p'}{p}}} (\disak)^{n+1}
\right)^{1-\frac{p'}{p}} .
\endmultline
$$
Now the exponent of $\disak$ in the last sum
is $\frac{\a' -\a}{\frac1{p'} - \frac1{p}} + n + 1 > n$ by
hypothesis, so we are done. Notice that if the sequence
$\sak_k$ was sparse, we could allow smaller values of $\a'$,
so that there is no hope of a general converse statement
analogous to the case (a).

The reasoning in the case $p=\infty$ is similar and
simpler.

(c) Our example will be a separated sequence which is as crowded
as possible (a net in the sense of \cite{Ro} or \cite{Lu}),
i.e. having the property that all points in the ball are less
than a constant invariant distance away from a point in the 
sequence.

Let $r \in (0,1)$. We pick a sequence 
$$
\{ a_{m,\ell} ,\ m \in \Bbb Z_+^*, 0 \leq \ell \leq L_m \},
$$
such that for any $\ell$, $a_{m,\ell} = 1 - r^m$,
and for each $m$ the set 
$\{ a_{m,\ell} , 0 \leq \ell \leq L_m \}$ is maximal in
the sphere of radius $1 - r^m$ for the property that the
invariant balls $E(a_{m,\ell} , r)$ be disjoint. It is easy
to check that this sequence is separated and that 
$L_m \approx r^{-nm}$.

We choose $v_{m,\ell}$ so that $u_m := |v_{m,\ell}|$ depends
only on $m$. Now for $0<p<\infty$,
$$
\sum_{m,\ell} (\disak)^{n+1+\a p} |v_{m,\ell}|^{p}
\approx \sum_m r^{-nm + m(n+1+\a p)} u_m^{p} ,
$$
and 
$$
\sum_{m,\ell} (\disak)^{n+1+\a' p'} |v_{m,\ell}|^{p'}
\approx \sum_m r^{ m(1+\a' p')} u_m^{p'}
\geq  \sum_m \left( r^{ m(1+\a p)} u_m^{p} \right)^{\frac{p'}{p}},
$$
since $\a' p' \leq \a p' + \frac{p'}{p} -1$. By choosing
$u_m$ appropriately, we can make the last sum diverge,
while $\sum_m r^{ m(1+\a p)} u_m^{p} < \infty$.

In the case $p=\infty$, simply taking $u_m = r^{-m \a}$,
we see that 
$\sum_m r^{ m(1+\a' p')} u_m^{p'} \geq \sum_m 1 =\infty$. \qed

\enddemo

\par
\ {}
\par

\centerline{\it \S 3. Inclusions for $Int(\bap)$.}

In this section we show that the inclusions given in Lemma 2.1 and
Lemma 2.2 are also verified by the corresponding spaces of interpolating
sequences. First we recall some known facts about Carleson measures
which will also be used in \S 4. 

For any $t>0$ and $\zeta_0\in\bo$ consider the {\it Carleson window}
with centre $\zeta_0$ and radius $t$ defined by
$C_t(\zeta_0)=\{z\in\bo : |1-\bar\zeta_0 z|<t\}$. A Borel measure
$\nu$ in $\bo$ is a {\it $q$-Carleson measure} if
$$
\nu(C_t(\zeta_0))=O(t^q)\quad\forall
t>0\quad\forall\zeta_0\in\partial\bo\ .
$$
A $n$-Carleson measure is simply called a {\it Carleson measure}.
What we call $q$-Carleson measures were studied in \cite{Am-Bo} where
they were called "Carleson measures of order $q/n$".

One of the main features of Carleson measures is the following.

\proclaim{\smc Lemma 3.1 (\cite{Hr} \cite{Ci-Wo})}
Let $q\geq p>0$ and $\a\geq-1/p$. Let $\mu$ be a positive measure. Then the
following are equivalent:
\roster
\item"(a)" $\left(\ibo |f(z)|^q d\mu(z)\right)^{\frac 1q}\leq c\Vert 
f\Vert_{p,\a}\quad\quad\forall f\in\bap$.
\item"(b)" $\mu(C_t(\zeta))=O(t^{(\frac{n+1}p+\a)q})\quad\forall
t>0\quad\forall \zeta\in\partial\bo$. 
\endroster
\endproclaim

We will be mainly interested in Carleson conditions for the measures $\sum_k
(\disak)^q\delta_{a_k}$, which have important relationships with the values
$$
K(\sak,p,q):=\sup_{k\in\Bbb Z_+}\sumjk \frac{(\disak)^p(\disaj)^q}
{\akbaj^{p+q}}\qquad p,q>0\ .
$$

\proclaim {\smc Lemma 3.2}
\roster
\item"(a)" If $\sum_k (\disak)^n \delta_k$ is a Carleson
measure, then for any $p>0$, 
$K(\{ a_k \} , p , n) < \infty$. 
\item"(b)" A positive
measure $\mu$ on $\bo$ is a Carleson measure if and only if there
exists some $\beta>n/2$ such that
$$
\sup_{b \in \bo } \ibo 
\frac{(1-|b|^2)^{2\beta - n}}
{|1-\bar b z|^{2\beta} } 
 d \mu (z)< \infty .
$$
\endroster
\endproclaim

(a) is an immediate consequence of \cite{Ma, Lemma
1.4}. (b) is a well-known result, which can be found in \cite{Ma, Lemma
1.2} and ultimately goes back to \cite{Ga, pag. 239}.

We now come to the main result of this section:

\proclaim{\smc Theorem 3.3}
Let $p,p'>0$ and $\a\geq -1/p$ , $\a'\geq -1/p'$, 
satisfying one of the following
conditions:
\roster
\item"(a)" $p\leq p'$ and $\nap<\napp$.
\item"(b)" $p\geq p'$ and $\a+1/p<\a'+1/p'$.
\endroster
Then $Int(\bap)\subset Int(B_{\a'}^{p'})$.

\endproclaim

In the special case where $n=1$, 
Seip \cite{Se} 
has proved a stronger result than Theorem 3.3 (b),
namely $Int ( B_\a^\infty ) = Int (B_{\a-\frac1p}^p )$
(see Appendix).
This suggests that the inequality $\a+1/p<\a'+1/p'$
in (b) is critical. 

In fact, if we take a sequence $a \subset \Bbb B^1 \times \{0\} 
\subset \Bbb B^n$, then we see by Lemma 1.8,
that $a \in Int  ( B_{\a}^p (\Bbb B^n))$ iff
$a \in Int  ( B_{\a+\frac{n-1}p}^p (\Bbb B^1))$, which, by
Seip's result, is the same as 
$Int  ( B_{\a +\frac{n}p}^\infty (\Bbb B^1))$. Since
we know from \cite{Se} that  $Int  ( B_\beta^\infty (\Bbb B^1))$
$\subsetneqq Int  ( B_{\beta'}^\infty (\Bbb B^1))$ 
when $\beta
< \beta'$, this shows
that $Int (\bap (\bo)) \neq Int (\bapp (\bo))$ when
$\nap \neq \napp$, and that $Int (\bap (\bo)) \not\subset Int (\bapp (\bo))$
when $\nap > \napp$. 
This shows that the inclusions
 in parts (a) and (b) of Theorem 3.3 are strict. 

\demo{\sl Proof}
Let $\sak_k$ be $\bap$-interpolating and take  $f_k \in \bap$ with
$(\disaj)^{\naf}f_k(a_j)$ $=\delta_{jk}$ and $\|f_k\|_{p,\a}\leq c$, for
some constant $c>0$ independent of $k$. Given $m>0$ define
$G_k(z)=g_k(z)\cdot f_k(z)$, where
$$
g_k(z)=\frac{(\disak)^{\naf+m}}{(1-\bar a_k z)^{\nafp+m} }\ .
$$
For
a given $\{ \lambda_k \}_k \in \ell^{p'}$, let $G:= \sum \lambda_k G_k$.
>From this definition it follows immediately that
$(\disak)^{\nafp} G (a_k) = \lambda_k$, and we need
to prove that $G \in B_{\a'}^{p'}$.

Assume first $\a>-1/p$.

Case $0<p'\leq 1$. Since $\| \cdot \|_{p',\a'}^{p'}$
satisfies the triangle inequality, it's enough to show:

{\sl Claim}: There exists $c>0$ independent of $k$ such that
$\|G_k\|_{p',\a'}\leq c$ for all $k$. 

{\sl Proof:} By definition of $G_k$, $\|G_k\|_{p',\a'}^{p'}$ equals
$$
I_k := 
\ibo\frac{(\disak)^{(\naf+m)p'}(\diz)^{(\a'-\a)p'}}{\akbz^{(\nafp+m)p'}}(\diz)^{\a
p'}|f_k(z)|^{p'} dm(z)\ . 
$$

 Case (a). Since $p\leq p'$, estimate $|f_k(z)|^{p'-p}$ by Lemma 1.3;  we
see that
$$
\multline
I_k \preceq
\|f_k\|_{p,\a}^{p'-p} \times \\
\times
\ibo\frac{(\disak)^{(\naf+m)p'}(\diz)^{\a' p'+(p-p')(\naf)-\a
p}}{\akbz^{(\nafp+m)p'}}(\diz)^{\a p} |f_k(z)|^{p} dm(z)\ ,
\endmultline
$$
which, since $(\naf+m)p'+[\a' p'+(p-p')(\naf)-\a
p]=(\nafp+m)p' $, shows that $\|G_k\|_{p',\a'}^{p'}\preceq
\|f_k\|_{p,\a}^{p'}$. 

Case (b). We may assume $p>p'$, and apply H\"older's 
inequality with $P=p/p'$ and $Q=\frac p{p-p'}$. 
Since by hypothesis $(\a'-\a)\pp>-1$, one has, for $m$
large enough,
$$
\multline
I_k \preceq
\left[\ibo\frac{(\disak)^{(\naf+m)\pp}
(\diz)^{(\a'-\a)\pp}}{\akbz^{(\nafp+m)\pp}}
dm(z)\right]^{\frac{p-p'}p} \|f_k\|_{p,\a}^{p'}\preceq\\
\preceq(\disak)^{[(\naf+m)\pp-(\nafp+m)\pp+(\a'-\a)\pp+n+1]\frac{p-p'}p}
\|f_k\|_{p,\a}^{p'}\approx\|f_k\|_{p,\a}^{p'} . 
\qed
\endmultline
$$

Case $p'>1$. First we give a useful estimate:
\proclaim{\smc Lemma 3.4}
Let $\sak_k$ and $\{g_k\}_k$ be as above, and let $A$ be such that
$(\frac{n+1}p+\a+m)A-n-1>-1$. Then
$$
\sum_k|g_k(z)|^A\preceq (\diz)^{A[(\naf)-(\nafp)]}\ .
$$
\endproclaim

This is a consequence of \cite{Ma}; more 
precisely, it follows from Lemma 4.1(d) below
with exponents $P:= A[(\nafp)-(\naf)]$ and
$Q:= A(\naf + m)$.


Case (a). Using in succession H\"older's inequality 
(with $1/p'+ 1/q'=1$), Lemma 3.4,
and Lemma 1.3 for $|f_k(z)|$, we obtain:
$$
\multline
\| G \|_{p',\a'}^{p'} =
\ibo\left|\sum_k\lambda_k g_k(z) f_k(z)\right|^{p'} (\diz)^{\a' p'}
dm(z)\preceq\\
\preceq \ibo \left[\sum_k
|g_k(z)|^{q'}\right]^{\frac{p'}{q'}}\left[\sum_k\mlk^{p'}|f_k(z)|^{p'}\right]
(\diz)^{\a' p'} dm(z)\preceq\\
\int\limits_{\bo}\sum_k \mlk^{p'} |f_k(z)|^p (\diz)^{\a p}
(\diz)^{p'(\naf-\frac{n+1}{p'}-\a' )+(p-p')(\naf)+\a' p' -\a p} dm(z) ,
\endmultline
$$
which is controlled by $\sum_k\mlk^{p'}$, since
$p'[(\naf)-(\nafp)]+(p-p')(\naf)+\a' p' -\a p=0$. 

Case (b). Let $1/p+1/q=1$. We first estimate 
$$
|G(z)|=\left|\sum_k\lambda_k g_k(z) f_k(z)\right|\leq\left(\sum_k
\mlk^{(1-\frac {p'}{p})q}|g_k(z)|^q\right)^{\frac
1q}\left(\sum_k\mlk^{p'}|f_k(z)|^p\right)^{\frac 1p}\ . 
$$ 
Then,
applying again H\"older's inequality with $P=p/p'$ and $Q=\frac p{p-p'}$
it follows that $\|G_k\|_{p',\a'}^{p'}$ is bounded by 
$$ 
\multline
\eightpoint \ibo\left[\sum_k \mlk^{(1-\frac {p'}p)q} |g_k(z)|^q
(\diz)^{q(\a'-\a)}\right]^{\frac{p'}q}\left[\sum_k \mlk^{p'} |f_k(z)|^p
(\diz)^{\a p}\right]^{\frac{p'}p} dm(z)\preceq\\ 
\preceq
\left[\ibo\left[\sum_k \mlk^{\frac {p-p'}p q}|g_k(z)|^q
(\diz)^{q(\a'-\a)}\right]^{\frac{p'p}{q(p-p')}}
dm(z)\right]^{\frac{p-p'}p}\times\\ \times \left[\ibo \sum_k \mlk^{p'}
|f_k(z)|^p (\diz)^{\a p} dm(z)\right]^{\frac{p'}p}\ . 
\endmultline 
$$
The second factor is controlled by $\sum_k
\mlk^{p'}\|f_k\|^{p'}_{p,\a}\preceq\sum_k\mlk^{p'}$. 
Taking $\gamma, \delta\geq 1$ with $1/\gamma+
1/\delta=1$
and applying H\"older's inequality with exponents 
$a=\frac{p'(p-1)}{p-p'}$ and $b=\frac{p'(p-1)}{p(p'-1)}$, 
and then Lemma 3.4, we can
bound the integral appearing in the first factor by
$$
\multline
\ibo \left[\sum_k |g_k(z)|^{\frac{qb}{\delta}}\right]^{\frac
ab}\left[\sum_k\mlk^{p'}|g_k(z)|^{\frac{qa}{\gamma}}(\diz)^{qa(\a'-\a)}\right]
dm(z) \preceq \\
\preceq
\ibo\sum_k\mlk^{p'}
(\disak)^{(\naf+m)\frac{qa}{\gamma}}\frac{(\diz)^{qa(\a'-\a)+\frac{qa}{\delta}(\naf-\frac{n+1}{p'}-\a')}}{\akbz^{(\nafp+m)\frac{qa}{\gamma}}}
dm(z)\ . 
\endmultline
$$
Since by hypothesis $qa(\a'-\a)=\pp(\a'-\a)>-1$, we can choose
$\delta>1$ so that $
qa(\a'-\a)+\frac{qa}{\delta}(\naf-\frac{n+1}{p'}-\a')>-1$ and the
integral is finite. Then, once more by Lemma 1.1, we see that
the integral is bounded by
$\sum_k\mlk^{p'}$. This concludes the case $\a>-1/p$.

We now turn to the case $\a = -1/p$. First we handle
the special situation $p'=p$:
\proclaim{\smc Lemma 3.5}
For any $\a' > -1/p$, $Int (H^p) \subset Int (B_{\a'}^p)$.
\endproclaim

Accepting this, suppose $(p',\a')$ satisfy
(a) in Theorem 3.3; then there exists $\a_1 > -1/p$
such that $\frac{n+1}p + \a_1 < \napp$. Likewise if $(p',\a')$ satisfy
(b) then there exists $\a_1 > -1/p$
such that $\frac{1}p + \a_1 < \uaprime$. In either
case, applying Lemma 3.5, then the case $\a > -1/p$
of Theorem 3.3, we get
$Int (H^p) \subset Int (B_{\a_1}^p) \subset Int (B_{\a'}^{p'})$.
\qed
\enddemo

\demo{\sl Proof of Lemma 3.5}
Define $f_k$, $g_k$, $G_k$ and $G$ as before for
$\{ \lambda_k \}_k \in \ell^p$.

Case $p \leq 1$. 
It is enough to prove that 
$\|G_k\|_{p,{\a'}}^p\leq c$ for all $k$, i.e.
$$
\ibo (\diz)^{\a' p} | g_k (z) |^p
|f_k (z) |^p d m (z) \leq c\quad \text{ for all } k.
$$
Applying Lemma 3.1 with $p=q$ and $\a=-1/p$, we see that it 
suffices to prove that, for an appropriate choice of the parameter $m$
in $g_k$,
$$
(\diz)^{\a' p}  | g_k (z) |^p d m (z)
$$
is a Carleson measure with Carleson norm independent of
$a_k$. 
To see this, we 
apply Lemma 3.2(b) with $2 \beta = n+1+\a' p +mp$. 
By the hypotheses on $\a'$ and $m$ we
have $2\beta > n+mp >n $. Then, by Lemma 1.1(c),
$$\align
\sup_{a,b \in \bo } & \ibo 
\frac{(1-|b|^2)^{2\beta - n}(\diz)^{\a' p} (1-|a|^2)^{mp + n}}
{|1-\bar b z|^{2\beta} |1-\bar a z|^{n+1+\a' p +mp}} 
 d m (z)\preceq\\
&\preceq
\sup_{a,b \in \bo } \,
\frac{(1-|b|^2)^{1+\a' p +mp} (1-|a|^2)^{mp + n} }
{|1-\bar a b|^{n+1+\a' p +mp}} \,
\{ \min( 1-|a|^2, 1-|b|^2 ) \}^{-mp} \ ,
\endalign
$$
which is finite since $\max ( 1-|a|^2, 1-|b|^2 ) \preceq |1-\bar a b|$.

Case $1<p\leq 2$.
By H\"older's inequality,
$$
\ibo (\diz)^{\a' p} |G(z)|^p  dm(z)\leq
\sum_j |\lambda_j |^p \ibo (\diz)^{\a' p} 
\left( \sum_k |g_k (z)|^q \right)^{\frac{p}{q}} | f_j (z)|^p dm(z) ,
$$
where $\frac1p + \frac1q = 1$.
Again by Lemmas 3.1 and 3.2(b), it's enough to consider
$$
\sup_{b \in \bo } \ibo 
\frac{(1-|b|^2)^{2\beta - n}}{|1-\bar b z|^{2\beta} } 
(\diz)^{\a' p} 
\left( \sum_k |g_k (z)|^q \right)^{\frac{p}{q}}  \, dm(z) \ . \tag{1}
$$
Since $\frac{p}{q}=p-1\leq1$, this integral is bounded by
$$
S := \sum_k \ibo 
\frac{(1-|b|^2)^{2\beta - n} (\diz)^{\a' p} (1-|a_k|^2)^{mp + n}}
{|1-\bar b z|^{2\beta} |1-\bar a_k z|^{n+1+\a' p +mp} } 
\, dm(z) \ .
$$
Choosing again $2 \beta = n+1+\a' p +mp > n+mp >n $ and
applying Lemma 1.1(c), we get $S\leq S_1 + S_2$, where
$$\align
S_1 &:= \sum_{k: \disak \leq 1-|b|^2}
\frac{(1-|b|^2)^{1+\a' p +mp} (1-|a_k|^2)^{n}}
{|1-\bar a_k b |^{n+1+\a' p +mp}  }  \\
S_2 &:= \sum_{k: \disak > 1-|b|^2}
\frac{(1-|b|^2)^{1+\a' p} (1-|a_k|^2)^{mp+n}}
{|1-\bar a_k b |^{n+1+\a' p +mp}  }  \ ,
\endalign
$$
so that $\sup_{b \in \bo} S_1 \leq K( \sak , 1+\a' p +mp , mp+n)$,
$\sup_{b \in \bo} S_2 \preceq K( \sak , 1+\a' p , n)$. Since
$\sak \in Int (H^p)$, $p>1$,
we know from \cite{Th1, Theorem 2.2} that $\sum_k (\disak)^n \delta_{a_k}$
is an $n$-Carleson measure, so Lemma 3.2(a) allows
us to conclude that both quantities are finite.

Case $2<p<\infty $. As before, it is enough to consider (1). 
We apply first H\"older's inequality with exponents
$\frac{p-1}{p-2}$ and $p-1$, then Lemma 3.4 to get
$$ 
\left( \sum_k |g_k (z)|^q \right)^{p-1}\leq 
\left( \sum_k |g_k (z)|^{\frac{p}{2(p-2)}} \right)^{p-2}
\sum_k |g_k (z)|^{\frac{p}2} 
\preceq
(\diz)^{-\frac{\a' p +1}2}
\sum_k |g_k (z)|^{\frac{p}2} \ ,
$$
so that in this case 
$$
S \preceq
\sup_{b \in \bo} \sum_k \ibo 
\frac{(1-|b|^2)^{2\beta - n} (\diz)^{\frac{\a' p -1}2} 
(1-|a_k|^2)^{\frac12 (mp + n)}}
{|1-\bar b z|^{2\beta} |1-\bar a_k z|^{\frac12 (n+1+\a' p +mp)} }
\, dm(z) \ .
$$
This time choose $2\beta = \frac12 (n+1+\a' p +mp)
> \max (\frac12 (mp + n) , n)$, which requires 
$m> \frac{n}p - (\a'+\frac1p )$. As above, 
$ S\leq S_1 + S_2$, where 
$$
\align
S_1 &:= \sum_{k: \disak \leq 1-|b|^2}
\frac{(1-|b|^2)^{\frac12 (-n+1+\a' p +mp)}(1-|a_k|^2)^{n}}
{|1-\bar a_k b |^{\frac12 (n+1+\a' p +mp)} }   \\
S_2 &:= \sum_{k: \disak > 1-|b|^2}
\frac{(1-|b|^2)^{\frac12 (1+\a' p)} (1-|a_k|^2)^{\frac12 (mp + n)} }
{|1-\bar a_k b |^{\frac12 (n+1+\a' p +mp)} } , \\
\endalign
$$
and requiring finally $mp \geq n$, we conclude as before.

Case $p=\infty$. We can estimate $\sup_{z\in\bo}(\diz)^{\a'}|G(z)|$
by a straightforward application of Lemma 3.4.
\qed
\enddemo

{\bf Remark.} We have actually shown a slightly stronger
result, namely that if functions $f_k$ exist with the properties
mentioned at the very beginning of the proof (in
the case where $\a =0$, $p=\infty$, this hypothesis
is sometimes called "uniform separatedness"), then $\sak \in Int (\bapp)$
for $(\a', p')$ verifying (a) or (b).

Notice that the proof of Theorem 3.3 cannot be used to show the
intuitive conjecture that $Int(H^p)\subset Int(H^{p'})$, for $p'<p$.
It is also interesting to note that the proof only uses
that $\sum_k (\disak)^n \delta_{a_k}$
is an $n$-Carleson measure in the case
$\sak \in Int (H^p)$, $p>1$, which follows
the arguments of \cite{Ca-Ga} and is much easier
to prove than for $p=1$ (see \cite{Th2, sec. 2.2}, while the case $p<1$ is not
known to us for $n\geq 2$.

\par
\ {}
\par

\centerline{\it \S 4. Sufficient conditions.}

In this last paragraph we give sufficient conditions for a sequence
$\sak_k$ to be $\bap$-interpolating in terms of the values $K(\sak,
p,q)$ defined in the previous section.

\proclaim{\smc Lemma 4.1 (\cite{Ma})}
The following conditions are equivalent:
\roster
\item "(a)" $\sak_k$ is the union of a finite number of separated sequences.
\item"(b)" $\ K(\sak,p,q)<+\infty \quad\forall q> n\ \forall p\leq q $.
\item "(c)" $\exists q\geq n,p \quad :\quad K(\sak,p,q)<+\infty $.
\item "(d)" $ \ \forall p>0\quad \forall
q>n\quad\sup_{z\in\bo}\sumk\frac{(\diz)^p(\disak)^q}{\akbz^{p+q}}<+\infty
$. 
\item "(e)" $\sum_k(\disak)^q\delta_{a_k}$ is a $q$-Carleson measure,
for $q>n$. 
\endroster
\endproclaim

We will consider, given a
sequence $\sak_k$, the restriction operator $T(f)=\{f(a_k)\}_k$ defined
in Lemma 1.5. Notice that
$$
\|T(f)\|_{\ell_{\nap}}^p=\sum_k |f(a_k)|^p(\disak)^{n+1+\a
p}=\int_{\bo}|f(z)|^p d\mu(z),
$$
where $\mu=\sum_k(\disak)^{n+1+\a p}\delta_{a_k}$. From Lemmas 3.1
and 4.1 we deduce thus that $T$ maps $\bap$ boundedly on
$\ell_{\nap}(\sak)$ if and only if $\sak_k$ is a finite union of
separated sequences. This gives a partial converse to Lemma 1.5.

The first result we give in this section deals with the case $p=1$.

\proclaim{{\smc Theorem 4.2}} Let $\sak_k$ be a separated sequence in
$\bo$. If there exists $m>0$ such that
$$
 K(\sak,m,n+1+\a)<1
$$ 
then $\sak_k$ is $B_\a^1$-interpolating.
\endproclaim
 
\demo{\sl Proof}
Consider $\
T:B_{\a}^{1}\longrightarrow\ell_{n+1+\a}^1(\sak)$ defined above.
In order to show that
$T$ is onto define, given $v=\{v_k\}\in\ell_{n+1+\a}^1(\sak)$, the
"approximate extension" 
$$
E(v)(z)=\sumk v_k\left(\frac{\disak}{1-\bar a_k z}\right)^{n+1+\a+m}\ .
$$
Using Lemma 1.1 it is immediately verified that $E(v)$ is 
in $B_{\a}^{1}$:
$$
\multline
\Vert E(v)\Vert_{1,\a}\leq \ibo(\diz)^\a\sumk
|v_k|\frac{(\disak)^{n+1+\a+m}}{\akbz^{n+1+\a+m}}\; dm(z)\leq\\
\leq\sumk|v_k|(\disak)^{n+1+\a+m}\ibo\frac{(\diz)^{\a}}{\akbz^{n+1+\a+m}}\;
dm(z)\approx\Vert v\Vert\ .
\endmultline
$$
On the other hand $TE-Id$, 
regarded as operator on $\ell_{n+1+\a}^1(\sak)$, has norm strictly
smaller than 1:
$$
\multline
\Vert TE(v)-v\Vert=\sumk(\disak)^{n+1+\a}|(TE(v))_k-v_k|\leq\\
\allowdisplaybreak
\leq\sumk(\disak)^{n+1+\a}\sumjk |v_j|\left(\frac{\disaj}{\akbaj}\right)^{n+1+\a+m}=\\
=\sumj(\disaj)^{n+1+\a}|v_j|\left(\sum_{k:k\neq j}\frac{(\disaj)^m
(\disak)^{n+1+\a}}{\akbaj^{n+1+\a+m}}\right)< 
\Vert v\Vert .
\endmultline
$$
Hence the series
$$
Id+\sumk (TE-Id)^k
$$
converges and defines an inverse to $TE$. The operator $\tilde
E=E(TE)^{-1}$ provides finally the inverse of $T$. \qed
\enddemo

Notice that by the invariance under automorphisms 
of the $\bap$-interpolating sequences the hypothesis in Theorem 4.2 
can be replaced by the seemingly weaker assumption of the existence 
of an automorphism $\phi_\zeta$ such that $K(\{\phi_\zeta(a_k)\},m,n+1+\a)<1$.

\proclaim{{\smc Corollary 4.3}}
\roster 
\item"(a)" Let $\sak_k$  be a separated sequence. 
There exists $\a>0$ such that $\sak_k$ is $B_{\a}^{1}$-interpolating.
\item "(b)" Let $\a>-1$. There exists $\delta\in (0,1)$ such that any
sequence $\sak_k$ verifying $d_G(a_j,a_k)\geq\delta$ for any $j\neq
k$ is $B_\a^1$-interpolating
\endroster
\endproclaim

\demo{\sl Proof} Lemma 4.1 shows that 
$$
K(\sak,n+\frac{\a+1}2,n+\frac{\a+1}2)
=\sup_{k}\sumjk\left(1-d_G^2(a_j,a_k)\right)^{n+\frac{\a+1}2}<+\infty
\ . 
$$ 
If $\delta$ is such that $
1-d_G^2(a_j,a_k)<1-\delta^2 $ we have then:
$$
\sup_{k\in\Bbb Z_+}\sumjk\left(1-d_G^2(a_j,a_k)\right)^{n+1+\a}\leq 
(1-\delta^2)^{\frac{\a+1}2} K(\sak,n+\frac{\a+1}2,n+\frac{\a+1}2)\ .
$$
In both cases (a) and (b) we can finish the proof by choosing respectively
$\a$ or $\delta$ so that $(1-\delta^2)^{\frac{\a+1}2} 
K(\sak,n+\frac{\a+1}2,n+\frac{\a+1}2)<1 $ and applying Theorem 4.2. \qed
\enddemo

Theorem 4.2 along with the following Mills' lemma provides also another
characterization of the sequences appearing in Lemma 4.1. 

\proclaim{{\smc Mills' lemma}}(See \cite{Ga} or \cite{Th1}).
Let $A_{jk}$, $\ j,k\in\Bbb Z_+$, be non-negative real numbers 
such that $A_{jk}=A_{kj}$ and $A_{jj}=0$ for any $j$ and $k$. 
If $\ \sup_{k\in\Bbb Z_+}\sum_{j\in\Bbb Z_+} A_{jk}=M<+\infty\ $, 
there exists a partition $\Bbb Z_+=S_1\cup S_2$, $\ S_1\cap S_2=
\emptyset$ such that
$$
\sup_{k\in S_i}\sum_{j\in S_i} A_{jk}\leq\frac M2\qquad i=1,2\ .
$$
\endproclaim

\proclaim{\smc Corollary 4.4} A sequence $\sak_k$ is the union of a 
finite number of separated sequences if and only if it is the union of
a finite number of $B_\a^1$-interpolating sequences.
\endproclaim

\demo{\sl Proof} The reverse implication is given directly by Corollary
1.6. To see the direct one we apply (b) of Lemma 4.1 with $p=q=n+1+\a$
and Mills' lemma with
$$
A_{jk}=\frac{(\disak)^{n+1+\a}(\disaj)^{n+1+\a}}{\akbaj^{2(n+1+\a)}}\
.
$$
For any $N\in\Bbb Z_+$ one can split $\sak_k$ into $2^N$ sequences
$\{b_k^l\}_k$, $l=1,...,2^N$, such that $K(\{b_k^l\},n+1+\a,n+1+\a)<
\frac{1}{2^N}K(\sak,n+1+\a,n+1+\a)$. Taking $N$ sufficiently large this
term becomes smaller than 1, what by Theorem 4.2 yields the stated
result. \qed

\enddemo

With the same methods it is also possible to obtain sufficient
conditions for a sequence to be $\bap$-interpolating, $p>1$. However,
these conditions are not so well adapted to the nature of the $\bap$
spaces, in the sense that they are symmetrical in $p$ and the
conjugated exponent $q$. In the proof, which goes like
\cite{Th1, Proposition 3.2}, we will use the duality between $\bap$ spaces.
Consider the product given by
$$
<f,g>=:\ibo f(z)\; \overline{g(z)}\; (\diz)^{\a p} dm(z)\ .
$$
Using Lemma 1.1 and some standard results for classical Bergman spaces
(see \cite{Am, Lemme 1.2.3} and \cite{Ru, chap. 7}) it is easy to
prove that in case $ 1<p<\infty$, the dual space of $\bap$ with respect to this
product is $B_{\beta}^{q}$, where $\frac 1p+\frac 1{q}=1$ and $\beta q=\a
p$. Furthermore, there is a
reproducing kernel for $\bap$ functions, namely
$$
K_z(\zeta)=\frac{\Gamma(n+\a p+1)}{\Gamma(n+1)\Gamma(\a p+1)}\quad
\frac{1}{(1-\zeta\bar z)^{n+1+\a p}}\ . 
$$

\proclaim{\smc Theorem 4.5} Let $1<p<\infty$ and let $q$ be its conjugated 
exponent. Let $\beta=\a p/q$. If there exist $c_1,c_2>0$ such that 
$c_1c_2<1$ and
\roster
\item " " $ K\bigl(\sak,\frac{n+1}p+\a,\frac{n+1}q+\beta\bigr)\leq
c_1^q\quad\quad $; 
$ \quad\quad
K \bigl(\sak,\frac{n+1}q+\beta,\frac{n+1}p+\a\bigl)\leq c_2^p $, 
\endroster
then $\sak_k$ is $\bap$-interpolating.
\endproclaim

\demo{\sl Proof}
Given $\bak\in\lap(\sak)$ take the approximate extension
$$
E(v)(z)=\sumk v_k\left(\frac{\disak}{1-\bar a_k z}\right)^{n+1+\a p}\ 
.
$$
Using the duality described above, the reproducing kernel for 
$\bap$ and Lemma 4.1 with
$\mu=\sum_k (\disak)^{n+1+\beta q}\delta_{a_k}$, one has that
$\|E(v)\|_{p,\a}\leq c\|v\|_{p,\nap}$.

On the other hand, if $T$ denotes the operator on $\bap$ associated to
$\sak_k$, we have $\Vert TE-Id\Vert<1$, since
$$
\multline
\sumk(\disak)^{n+1+\a p}|(TE(v)-v)_k|^p\leq\\
\allowdisplaybreak
\leq\sumk\left(\sumjk\frac{(\disak)^{\frac{n+1}p+\a}
(\disaj)^{\frac{n+1}q+\beta}}{\akbaj^{n+1+\a p}}\right)^{\frac pq}\times\\
\allowdisplaybreak
\times\sumjk\frac{(\disak)^{\frac{n+1}p+\a}(\disaj)^{\frac{n+1}q+\beta}}
{\akbaj^{n+1+\a p}}(\disaj)^{n+1+\a p}|v_j|^p\leq\\
\leq c_1^p\sumj(\disaj)^{n+1+\a p}|v_j|^p\sum_{k:k\neq 
j}\frac{(\disaj)^{\frac{n+1}q+\beta}(\disak)^{\frac{n+1}p+\a}}{\akbaj^{n+1+\a
p}}\leq(c_1c_2)^p\Vert v\Vert_{p,\nap}\ .
\endmultline
$$
This shows that $TE $ is invertible and, as in Theorem 4.2, $T$ is onto. \qed
\enddemo

Similar corollaries to 4.3 and 4.4 can be derived from Theorem 4.5,
with some restrictions on the values $p$ and $\a$. 

\proclaim{\smc Corollary 4.6} Let $p>1$, $q$ be the conjugated 
exponent and denote $A(\alpha,p)=(n+1+\alpha p)\min(1/p,1/q)$. 
\roster
\item"(a)"  If there exists $c_0<1$ such that
$$
K(\sak,A(\alpha,p),A(\alpha,p))\leq 2^{-(n+1+\a p)|1-\frac 2p|} c_0
$$
then $\sak_k$ is $\bap$-interpolating.
\item"(b)" Let $\sak_k$ be separated. There exists $\a>0$ such that 
$\sak$ is $\bap$-interpolating.
\item"(c)" Let $A(\a,p)>n$. There exists $\delta>0$ such that any
sequence $\sak_k$ verifying $d_G(a_j,a_k)\geq\delta$ for all $j\neq k$
is $\bap$-interpolating.
\item"(d)" Let $A(\a,p)>n$. Then
$\sak_k$ is a finite union
of separated sequences if and only if it is a finite union of 
$\bap$-interpolating sequences.
\endroster
\endproclaim
 
Part (c), like (b) of Corollary 4.3, is a particular case of a
theorem of Rochberg which actually shows that the result holds for
any $p>0$ and $\a>-1/p$ (see \cite{Ro, pag. 231}).

\demo{\sl Proof} (a) is the analog of \cite{Th1, Corollary 3.3}.
(b), (c) and (d) are derived from (a), like in Corollaries 4.3 and 4.4. \qed
\enddemo

In order to use Theorems 4.2 and 4.5 (and their corollaries) for a given 
sequence $\sak_k$ it can be useful to know whether it is enough to 
verify the conditions given therein for a sequence obtained from the 
original one by taking off a finite number of points. This is 
equivalent to asking whether the process of adding a finite number of 
points to a $\bap$-interpolating sequence gives again a sequence that 
is $\bap$-interpolating. 

\proclaim{\smc Theorem 4.7}
The union of a $\bap$-interpolating sequence and a finite number of 
points is again $\bap$-interpolating.
\endproclaim

\demo{\sl Proof} It is enough to show that the union of a
$\bap$-interpolating sequence and one point is $\bap$-interpolating, and
by invariance under automorphisms of the $\bap$-interpolating sequences,
we can assume that this point is 0. Let then $\sak_k$ be the original
sequence and let $\delta>0$ be such that $|a_k|\geq\delta$, for all $k$. 

We claim first that it is enough to find $f\in\bap$ with $f(a_k)=0$ for
all $k$ and $f(0)\neq 0$.  To see this let $\sbk\cup
v_0\in\ell_{\frac{n+1}p+\alpha}^p(\sak\cup 0)$, and let $g\in\bap$ be
such that $g(a_k)=v_k$. Then the function
$$
F(z)=g(z)+\frac{v_0 -g(0)}{f(0)}\; f(z)
$$
belongs to $\bap$ and  $F(0)=v_0$, $F(a_k)=v_k$ for all
$k$.

Suppose that all $f\in\bap$ with $f(a_k)=0$ for all $k$ have $f(0)=0$.
This implies that for any $f\in\bap$, the value $f(0)$ is determined by
the values $f(a_k)$, since the difference of two functions with the same
values on $\sak_k$ vanishes at 0. 

Assume $1\leq p <\infty$ and define the functional $\Lambda:
\ell^p\longrightarrow\Bbb C$ by
$$
\Lambda(\sbk)=f(0)\ ,
$$
where $f\in\bap$ is such that $f(a_k)=(\disak)^{-(\frac{n+1}p+\alpha)}
v_k$ for all $k$. Since $f(0)$ is determined only by these values,
 which are actually independent of $f$, we have that
 $\Lambda$ is linear. It is also continous: 
$$
|\Lambda(b)|=|f(0)|\leq c \Vert f\Vert_{p,\alpha}\leq c M\Vert
\{(\disak)^{-(\frac{n+1}p+\alpha)}
v_k\}\Vert_{p,\frac{n+1}p+\alpha}=c M\Vert v\Vert_p\ ,
$$
$M$ denoting the interpolation constant of $\sak_k$. So
$\Lambda\in\ell^q$, in the sense that there exists $\{c_k\}_k\in\ell^q$
such that
$$
\Lambda(v)=\sumk v_k  c_k\qquad\forall\ v=\sbk\in\ell^p\ .
$$
Consider now the sequences $v^j=\{\delta_{jk}\}_k\in\ell^p$ and a
function $f_j\in\bap$ with
$f_j(a_k)=(\disak)^{-(\frac{n+1}p+\alpha)}\delta_{jk}$. By definition
$$
\Lambda(v^j)=f_j(0)=\sumk\delta_{jk} c_k= c_j\ .
$$
Take now the functions $F_j(z)=(|a_j|^2-\bar a_j z)f_j(z)$. Obviously
$F_j\in\bap$ and $F_j(a_k)=0$ for all $k$. Therefore
$F_j(0)=|a_j|^2f_j(0)=|a_j|^2 c_j=0$, hence, $c_j=0$. This shows that
$\Lambda\equiv 0$, which is evidently false, since there are many
functions in $\bap$ not vanishing at 0. 

The case $0<p<1$ is solved in the same way, using that the dual of
$\ell^p$ is $\ell^\infty$. For the case $p=\infty$ we can restrict the
functional $\Lambda$ to the subspace $c_0\subset\ell^\infty$ of
sequences with limit 0 and apply the same argument. \qed

\enddemo

 \par
\ {}
\par

\centerline{\it \S 5. Appendix.}

Given a sequence $\sak_k\subset\Bbb B^1$ consider its upper uniform density
$$
D(\sak)=\limsup_{r\to 1} \sup_{z\in\Bbb B^1}
\frac{\sum_{1/2<|\phi_z(a_k)|<r}\log \frac 1{|\phi_z(a_k)|}}
{\log \frac 1{1-r}}\ .
$$

\proclaim{\smc Lemma 5.1} Let $\sak_k\subset \Bbb B^1$.
Then $\sak\subset Int(\bap)$ if and only if $D(\sak)<\a+1/p$.
\endproclaim

\demo{\sl Proof} Assume first that $D(\sak)<\a+1/p$ and define
$\epsilon=1/2(\a+1/p-D(\sak))$. There exist functions
$f_k\in B_{\a+1/p-\epsilon}^\infty$ and $C>0$ with 
$\|f_k\|_{\infty,\a+1/p-\epsilon}\leq C$ for all $k$ and 
$f_k(a_j)=\delta_{jk}(\disak)^{-(\a+1/p)+\epsilon}$
(see \cite{Se, pag 34}). As in Theorem 3.3 
(see the remark at the end of the proof) this implies that
$\sak_k$ is $\bapp$-interpolating, for any $(p',\a')$ such that
$\a'+1/p'>\a+1/p-\epsilon$, and in particular for $\a$ and $p$.

Assume now that $\sak_k$ is $\bap$-interpolating.

\proclaim{\smc Lemma 5.2} Let $z\in\Bbb B^1$ such that
$d(z,\sak)\geq\delta_0$. There exists $f\in\bap$ with $f(a_k)=0$ for
all $k$, $f(z)=1$ and  $\|f\|_{p,\a}\leq 1+ M\delta_0^{-1}$, where
$M$ denotes the constant of interpolation of $\sak$.
\endproclaim

\demo{\sl Proof} By invariance under automorphisms we can suppose
that $z=0$. There exists $f_0\in \bap$ such that $f_0(a_1)=1/a_1$, 
$f_0(a_j)=0$ for all $j\geq 2$ and $ \|f_0\|_{p,\a}\leq
M\delta_0^{-1}$. Then the function $f(z)=1-zf(z)$
satisfies all the requirements. \qed 
\enddemo

Since $M_p^p(f,r)=\int_0^{2\pi} |f(re^{i\theta})|^p d\theta$ is an
increasing function of $r$,
$$
\|f\|_{p,\a}^p =\int_0^1 2\pi M_p^p(f,r) (1-r^2)^{\a p} r\; dr\geq
 M_p^p(f,r_0) \frac {\pi}{\a p+1} (1-r_0^2)^{\a p+1}\ .
$$
Hence, $M_p^p(f,r)\preceq \|f\|_{p,\a}^p (1-r^2)^{-(\a p+1)}$. 
By Jensen's inequality
$$
\exp\left( p\int_0^{2\pi}\log |f(re^{i\theta})|\frac{d\theta}{2\pi}\right)
\leq M_p^p(f,r)\leq C (1-r^2)^{-(\a p+1)}\ ,
$$
and thus
$$
\int_0^{2\pi}\log |f(re^{i\theta})|\frac{d\theta}{2\pi}\leq
\frac 1p \log C +(\a +\frac 1p)\log\bigl(\frac
1{1-r^2}\bigr)\ . 
$$
>From Jensen's formula it
follows now that 
$$
\sum_{1/2<|\phi_z(a_k)|<r}\log \frac r{|\phi_z(a_k)|}
\leq \frac 1p
\log C +(\a +\frac 1p)\log\bigl(\frac 1{1-r^2}\bigr)\ ,
$$
and therefore $D^+(\sak)\leq \a +1/p$. To see that the inequality is strict
take a sequence $\{a_k'\}_k$ and $\delta_0$ so that $d(a_k,a_k')<
\delta_0$ for all $k$ and $(1+\delta_0) D(\sak)\leq D(\{a_k'\})$ (see
the proof of \cite{Se, Lemma 6.6}). An application of the argument
above to the sequence $\{a_k'\}_k$ shows finally that
$D(\sak)< \a +1/p$. \qed
\enddemo

\enddocument